\documentclass[10pt]{article} 
\usepackage{amssymb}
\usepackage{amsmath, amsfonts, amsthm, amssymb, verbatim,a4wide}  
\theoremstyle{plain}
        \newtheorem{theorem}{Theorem}[section]

        \newtheorem{proposition}[theorem]{Proposition}
        \newtheorem{lemma}[theorem]{Lemma}
        
        \newtheorem{remark}[theorem]{Remark}

\numberwithin{equation}{section}
\newcommand \Ubar {\widetilde U}
\newcommand \demi   {{1/2}}
\newcommand \alphab {\underline \alpha}
\newcommand \ZZ     {\mathbb{Z}}
\newcommand \be         {\begin{equation}}
\newcommand \ee         {\end{equation}} 
\newcommand \del        \partial
\newcommand \eps    \varepsilon

\newcommand \im {\text{Im }}    
    
\newcommand \RR   	{{\mathbb R}}

\newcommand \Dcal    {\mathcal D}

\newcommand \Rn    {\mathbb{R}^m} 

\newcommand \RN    {\mathbb{R}^N} 
\newcommand \Ocal   {\mathcal O}    
\newcommand \Ecal   {\mathcal E}    
   
\newcommand \Mcal {\mathcal M}
\newcommand \Lcal {\mathcal L}

\newcommand \Abar  {{\overline A}}

\newcommand \Rcal {\mathcal R}
\newcommand \lam \lambda

\let\oldmarginpar\marginpar
\renewcommand\marginpar[1]{\-\oldmarginpar[\raggedleft\footnotesize #1]%
{\raggedright\footnotesize #1}}

\begin{document} 
%
\title{A framework for late-time/stiff relaxation asymptotics} 
\author{Philippe G. LeFloch\footnote{Laboratoire Jacques-Louis Lions, Centre National de la Recherche Scientifique, 
Universit\'e Pierre et Marie Curie (Paris 6), 4 Place Jussieu, 
75252 Paris, France. 
Email : {\tt pglefloch@gmail.com.} 
\newline
Blog: {\tt philippelefloch.wordpress.com.} 
\newline To appear in: {\sl
Proc. " CFL condition - 80 years gone by'', Rio de Janeiro Aug. 2010, C.S. Kubrusly et al. (editors), 
Springer Proc. Math., 2011.} 
}
}  
\date{\today}  
\maketitle

\begin{abstract}  We consider solutions to nonlinear hyperbolic systems of balance laws
with stiff relaxation and formally derive a parabolic-type effective system describing the late-time asymptotics of these solutions. 
We show that many examples from continuous physics fall into our framework,
including the Euler equations with (possibly nonlinear) friction.  
We then turn our attention to the discretization of these stiff problems and introduce a new finite volume 
scheme which preserves the late-time asymptotic regime.  Importantly, our scheme requires only the classical 
CFL (Courant, Friedrichs, Lewy) condition associated with the hyperbolic system under consideration, rather than the more restrictive, 
parabolic-type stability condition. 
\end{abstract}

\tableofcontents
\newpage 


\section{Introduction} 

This short presentation is based on the joint work \cite{BLT} in collaboration with C. Berthon and R. Turpault. 
We are interested in hyperbolic models arising in continuum physics and, especially, describing 
complex multi-fluid flows involving several time-scales. The partial differential equations under consideration 
are nonlinear hyperbolic systems of balance laws with stiff relaxation sources. We investigate here
 the late-time behavior of entropy solutions. 

Precisely, our objective is to derive the relevant effective system ---which turns out to be of parabolic type--- 
and to investigate the role of a convex entropy associated with the given system of balance laws. 
As we show it, many examples from continuous physics fall into our framework. In addition, 
we investigate here the discretization of such problems, and we propose a new finite volume 
scheme which preserves the late-time asymptotic regime identified in the first part of this paper.  
 
An outline of this paper follows. In Section~2, we present a formal derivation of the 
effective equations associated with our problem. In Section~3, we demonstrate that many examples of continuous physics are covered by our theory. Finally, in Section~4, we are in a position to present the new discretization method and state several properties of interest. 


\section{Late-time/stiff relaxation framework}  

\subsection{Hyperbolic systems of balance laws} 

We consider systems of partial differential equations of the form 
\be
\label{100}
\eps \, \del_t U + \del_x F(U) = -{R(U) \over \eps}, \qquad  U=U(t,x) \in \Omega \subset \RN, 
\ee
where $t >0, \, x \in \RR$ denote the (time and space) independent variables. 
We make the following standard assumptions. 
The flux $F: \Omega \to \RN$ is defined on a convex and open subset $\Omega$. 
The first-order part of \eqref{100} is an hyperbolic system, that is, the matrix-valued field 
 $A(U) := D_UF(U)$ admits real eigenvalues    
and a basis of eigenvectors.

We are interested in the singular limit problem $\eps \to 0$ in the limit of late-time and stiff relaxation. 
In fact, two distinct regimes for systems like \eqref{100} can be considered. 
In the hyperbolic-to-hyperbolic regime, one replace $\eps \del_t U$ by $\del_t U$, and 
one establishes that solutions to 
$$ 
\del_t U + \del_x F(U) = - {R(U) \over \eps}, \qquad  U=U(t,x), 
$$
are driven by an effective system of equations ($\eps \to 0$) of {\sl hyperbolic type.} 
Such a study was pioneered by 
Chen, Levermore, and Liu \cite{CLL}. 
On the other hand, in the hyperbolic-to-parabolic regime under consideration in the present work, we 
obtain effective equations of {\sl parabolic type.} Earlier work by Marcati et al. \cite{Marcati1} discussed 
this regime too and established several important convergence theorems.

Our objective here is to introduce a general framework to deal with such problems. We make the following assumptions. 

\

\noindent{\bf  Assumption 1.} There exists an $n \times N$ matrix $Q$ with (maximal) rank $n<N$ such that
$$
Q R(U) = 0, \qquad U \in \Omega.
$$
Hence, $QU \in Q\Omega =: \omega $ satisfies
$$
\eps \, \del_t \big( QU  \big)  + \del_x \big( QF(U) \big) = 0. 
$$
 
\

\noindent{{\bf  Assumption 2.} There exists a map $\Ecal:\omega \subset \RR^m \to\Omega$ describing the equilibria $u \in \omega$, with   
$$
R(\Ecal(u))=0,
\qquad 
u = Q \, \Ecal(u).
$$
It is then convenient to introduce the equilibrium submanifold $\Mcal:= \big\{ U= \Ecal(u) \big\}$. 

\

\noindent{{\bf  Assumption 3.}  It is also assumed that 
$$
QF(\Ecal(u))= 0, \qquad  u\in\omega.
$$
To motivate this condition we observe that, formaly at least, the term 
$\del_x \big( QF(\Ecal(u)) \big)$ must vanish identitically, so that $QF(\Ecal(u))$ must be a constant, 
conveniently normalized to be $0$. 
 
\

\noindent{{\bf  Assumption 4.} For all $u\in \omega$, 
$$ 
\aligned
\dim\Big( \ker(B(\Ecal(u)))\Big)&=n,
\\
\ker\big( B(\Ecal(u)) \big) \cap \im\big( B(\Ecal(u))\big) & =\{0\}.
\endaligned
$$ 
Hence, the $N\times N$ matrix $B:=DR_U$ has ``maximal'' kernel on the equilibrium manifold.

\subsection{Chapman-Engskog-type expansion}

We proceed by using a Chapman-Engskog expansion in order to derive
effective equations satisfied by the local equilibria $u= u(t,x) \in \omega$. So, we write 
$$
U^\eps = \Ecal(u) + \eps \, U_1 +\eps^2 \, U_2 + \ldots, 
\qquad 
\quad 
 u:= QU^\eps,  
$$
which should satisfy  
$
\eps \, \del_t U^\eps + \del_x F(U^\eps) = - R(U^\eps)/\eps. 
$
It follows that 
$$
QU_1= QU_2 = \ldots = 0.
$$
For the flux we find 
$$
\begin{aligned}
F(U^\eps) = F(\Ecal(u)) & + \eps \, A(\Ecal(u))\,  U_1
+ 
\Ocal(\eps^2),
\end{aligned}
$$
and for the relaxation 
$$
\aligned
\frac{R(U^\eps)}{\eps}
= B(\Ecal(u))\,  U_1 & +
   \frac{\eps}{2}D^2_U R(\Ecal(u)).(U_1,U_1) 
 +
   \eps B(\Ecal(u))\,  U_2 + \Ocal(\eps^2).
\endaligned
$$
 In turn, we deduce that 
$$
\begin{aligned}
& \eps \, \del_t \big( \Ecal(u) \big)
+ \del_x \big(  F(\Ecal(u)) \big) + \eps \, \del_x \Big( A(\Ecal(u))\,  U_1 \Big) 
\\ 
&
= -B(\Ecal(u))\,  U_1 -\frac{\eps}{2} D^2_U R(\Ecal(u)).(U_1,U_1)
  -\eps B(\Ecal(u))\,  U_2 +\Ocal(\eps^2).
\end{aligned}
$$
 
We begin by considering the 
{\sl zero-order terms} and thus deduce that $U_1 \in \RN$ satisfies the linear system
$$
B(\Ecal(u))\,  U_1 = - \del_x \big( F(\Ecal(u))\big) \in \RN.
$$ 
We can solve this equation in $U_1$ provided we recall  that $QU_1 =0$ and observe the following. 
 
\begin{lemma}[Technical lemma] Let $C$ be an $N \times N$ matrix satisfying 
$\dim \ker C = n$, and $\ker C \cap \im C = \big\{ 0 \big\}$, and let $Q$ be an
$n \times N$ matrix of rank $n$. Then, for all $J\in\RR^N$, the system
$$
\begin{aligned}
& C\,  V =J,\\
& QV =0,
\end{aligned}
$$
admits a unique solution $V\in\RR^m$ if and only if $QJ=0$.
\end{lemma}

We can thus conclude with the following. 

\begin{proposition}[First-order corrector problem]
The first-order corrector $U_1$ is characterized uniquely by  
$$
\begin{aligned}
B(\Ecal(u))\,  U_1 & = -\del_x\big( F(\Ecal(u)) \big),
\\
QU_1 & = 0.
\end{aligned}
$$
\end{proposition}

We now turn our attention to the {\sl first-order terms} and we arrive at 
$$
\del_t \big( \Ecal(u) \big) + \del_x \Big( A(\Ecal(u))\,  U_1 \Big)
= -\frac{1}{2} D^2_U R(\Ecal(u)).(U_1,U_1)- B(\Ecal(u))\,  U_2.
$$
Multiplying by $Q$ and using $Q\Ecal(u)=u$ we find 
$$
\del_t u  + \del_x \Big( Q \, A(\Ecal(u))\,  U_1 \Big)
= -\frac{1}{2} Q \, D^2_U R(\Ecal(u)).(U_1,U_1) - Q \, B(\Ecal(u))\,  U_2.
$$
But, by differentiating the identity $Q R(U)=0$, it follows that 
$$
 Q \, D^2_U R . (U_1,U_1) \equiv 0, 
\qquad
 Q \, B \,  U_2 \equiv 0.
$$ 

\begin{theorem}[Late time/stiff relaxation effective equations] One has 
$$
\del_t u = -\del_x \Big( Q A(\Ecal(u))\, U_1 \Big) =: \del_x \left(\Mcal(u) \del_x u \right) 
$$
for some $n\times n$ matrix $\Mcal(u)$, 
where $U_1$ is the unique solution to the first-order corrector problem
$$
\begin{aligned}
B (\Ecal(u))\,  U_1 &= - A(\Ecal(u)) \del_x\big( \Ecal(u) \big),
\\
QU_1 & = 0.
\end{aligned}
$$
\end{theorem}


\subsection{The role of a mathematical entropy} 

Next, we investigate the consequences of assuming the existence of a mathematical entropy 
$\Phi:\Omega\to\RR$, satisfying by definition: 

\

\noindent{{\bf  Assumption 5.} There exists an entropy-flux $\Psi:\Omega\to\RR$ such that 
$$
D_U \Phi \,  A = D_U \Psi \qquad \text{ in } \Omega.  
$$
So, all smooth solutions satisfy
$$
\eps\del_t \big( \Phi(U^\eps) \big) + \del_x \big( \Psi(U^\eps) \big)
= 
- D_U \Phi(U^\eps)\,  {R(U^\eps) \over \eps}
$$
and, consequently, the matrix $D^2_U \Phi \, A$ is symmetric in $\Omega$.
In addition, we assume that the map $\Phi$ is convex, i.e.~the $N\times N$ 
matrix $D^2_U \Phi$ is positive definite on $\Mcal$.

\

\noindent{{\bf  Assumption 6.} The following compatibility property with the relaxation term holds:  
$$
\aligned
&  D_U \Phi \,  R \ge 0 \qquad  \text{ in } \Omega, 
\\
& D_U \Phi(U) = \nu(U) Q \in \RN, \qquad \nu(U) \in \Rn.  
\endaligned
$$

Returning to the effective equations 
$$
\del_t u= \del_x \Dcal,  
\qquad 
\Dcal := -Q  A(\Ecal(u))\,  U_1 
$$
and multiplying it by the Hessian of the entropy, we conclude that the term 
$U_1\in\RR^N$ is now characterized by 
$$
\aligned
\Lcal(u) U_1 &= - \big( D^2_U \Phi \big) (\Ecal(u)) \del_x \big( F(\Ecal(u)) \big),
\\
Q U_1 &= 0, 
\endaligned
$$
where $\Lcal(u) = D^2_U \Phi(\Ecal(u))B(\Ecal(u))$. 

Then, using the notation $\Lcal(u)^{-1}$ for the generalized inverse with constraint and setting
$$
S(u) := Q A(\Ecal(u)),
$$ 
 we obtain  
$$
\Dcal = S \Lcal^{-1} \big( D^2_U \Phi \big) (\Ecal) \del_x\big( F(\Ecal) \big).
$$
Finally, one can check that 
$$
\big( D^2_U \Phi \big) (\Ecal) \del_x \big( F(\Ecal) \big) = S^T  v,
$$
with 
$v := \del_x \left(  D_u \Phi(\Ecal) \right)^T $.

\begin{theorem}[Entropy structure of the effective system]
When the system of balance laws is endowed with a mathematical entropy, the effective equations take the form 
$$
\del_t u = \del_x \Big( L(u) \, \del_x \big( D_u \Phi(\Ecal(u)) \big)^T \Big), 
$$
where
$$
\aligned
L(u) & := S(u)\Lcal(u)^{-1}S(u)^T, 
\\
S(u) &:= Q A(\Ecal(u)), 
\\
\Lcal(u) &:= \big( D^2_U\Phi \big) (\Ecal(u)) B (\Ecal(u)),
\endaligned 
$$
where, for all $b$ satisfying $Qb =0$,  the unique solution to
$$
\Lcal(u) V =b, \quad \quad QV =0 
$$
is denoted by $\Lcal(u)^{-1} b$ (generalized inverse).  
\end{theorem}

Alternatively, the above result can be reformulated in terms of the so-called entropy variable 
$\big( D_u \Phi(\Ecal(u)) \big)^T$. 
Furthermore, an important dissipation property can be deduced from our assumptions, as follows. 
From the entropy property and the equilibrium property $R(\Ecal(u))=0$,  we find 
$$
\aligned
D_U \Phi R & \ge 0 \qquad  \text { in} \Omega,
\\
\left(D_U \Phi R \right)|_{ U=\Ecal(u)}& =0 \qquad  \text {in } \omega.
\endaligned
$$
Thus, the matrix $D^2_U\Big( D_U \Phi R \Big)|_{ U=\Ecal(u)}$ is non-negative definite. 
 It follows that 
$$ 
D^2_U\Big(D_U \Phi R \Big) 
 = D^2_U\Phi B + \left( D^2_U\Phi B \right)^T \qquad \text{ when } U=\Ecal(u),
$$
so that 
$$
D^2_U\Phi \, B |_{ U=\Ecal(u)} \geq 0 \qquad  \text {in } \omega. 
$$

The {\sl equilibrium entropy} $\Phi(\Ecal(u))$ has the property that its associated 
(entropy) flux $u \mapsto \Psi(\Ecal(u))$ is constant on the equilibrium manifold $\omega$.
Indeed, for the map $\Psi(\Ecal)$, we have 
$$
\aligned
D_u \big( \Psi(\Ecal) \big) & = D_U \Psi(\Ecal) D_u \Ecal
\\
&= D_U \Phi(\Ecal) A(\Ecal) D_u \Ecal.
\endaligned
$$
Observing that $\big( D_U \Phi \big) (\Ecal) = D_u \big( \Phi(\Ecal)\big) Q$, we obtain 
$$
\aligned
D_u \Big( \Psi(\Ecal(u)) \Big) &= D_u \Phi(\Ecal(u)) Q  A(\Ecal(u)) D_u
\Ecal(u)
\\
&=  D_u \Big( \Phi(\Ecal(u)) \Big) D_u QF(\Ecal(u)). 
\endaligned
$$
Since $QF(\Ecal) = 0$, then $D_u QF(\Ecal) =0$ and the proof is completed.

Therefore, $D_u \big( \Psi(\Ecal(u)) \big)=0$ for all $u \in \omega$.

Recalling the asymptotic expansion
$$
U^\eps = \Ecal(u) + \eps U_1 + ...,
$$
where $U_1$ is given by the first-order corrector problem, we write 
$$
\Psi(U^\eps)=\Psi(\Ecal(u)) +\eps \, D_U\Psi(\Ecal(u)) \, U_1 + \Ocal(\eps^2),
$$
and deduce  
$$
\del_x \Psi(U^\eps) = \eps \, \del_x D_U\Psi(\Ecal(u)) \, U_1 + \Ocal(\eps^2).
$$
Similarly, for the relaxation source term,  we have 
$$
D_U \Phi(U^\eps)R(U^\eps) = \eps^2D^2_U \Phi(\Ecal(u)) D_U
R(\Ecal(u)) U_1 + \Ocal(\eps^3).
$$
We thus obtain
$$
\aligned
& \del_t \big( \Phi(\Ecal(u)) \big)+ \del_x \Big( D_U \Psi(\Ecal(u)) \, U_1\Big) 
\\
& = - U_1^T \,  \left( D^2_U \Phi(\Ecal(u)) B(\Ecal(u)) \right) U_1. 
\endaligned
$$ 
But, we have already established 
$$
X\,  \, \big( D^2_U \Phi \big)(\Ecal) B(\Ecal) \, X \ge 0,
\qquad  X \in \RR^N.
$$ 

\begin{proposition}[Monotonicity of the entropy]
The entropy is non-increasing in the sense that 
$$
\del_t \big( \Phi(\Ecal(u)) \big) + \del_x \left( D_U \Psi(\Ecal(u)) \, U_1\right) \leq 0.
$$
\end{proposition}

In the notation given earlier, one thus have 
$$
\del_t \big( \Phi(\Ecal(u)) \big) = \del_x \Big( \big( D_u \big(\Phi(\Ecal(u) \big) \big) 
L(u) \, \del_x \big( D_u \big(\Phi(\Ecal(u) \big) \big)^T \Big).  
$$

\begin{remark} 
\label{rmk26}
The above framework can be extended to handle certain {\sl nonlinear diffusion regime,}
corresponding to the scaling 
$$
\eps \, \del_t U+\del_x F(U) = -\frac{ R(U)}{\eps^q}.
$$
The parameter $q \ge 1$ introduces an {\sl additional scale}, and is
indeed necessary for certain problems where the relaxation is nonlinear. 
The relaxation term is supposed to be such that 
$$
R\big(\Ecal(u)+\eps \, U\big) = \eps^q R\Big(\Ecal(u)+M(\eps) \, U\Big),
\qquad U\in\Omega, \quad u\in\omega, 
$$
for some matrix $M(\eps)$.  In that regime, the effective equations turn out to be of {\em nonlinear parabolic} type.   
\end{remark} 


\section{Examples from continuum physics} 

\subsection{Euler equations with friction term}
   
The simplest example of interest is provided by the Euler equations of compressible fluids with friction: 
\be
\label{Euler1} 
\begin{aligned}
& \eps \, \del_t\rho +\del_x( \rho v ) =0
\\
& \eps \, \del_t( \rho v) +\del_x \big( \rho v^2 + p(\rho) \big) = - {\rho v \over \eps}
\end{aligned}
\ee
in which the density $\rho \geq 0$ and the velocity component $v$ represent the main unknowns while
the pressure $p:\RR^+\to\RR^+$ is prescribed and satisfy the hyperbolicity condition $p'(\rho)>0$ for all $\rho>0$. 
Then, the first-order homogeneous system is strictly hyperbolic and \eqref{Euler1} fits into our
late-time/stiff relaxation framework provided we set  
$$
U=\left(\begin{array}{c}
\rho \\ \rho v
\end{array}\right),
\quad
F(U)=\left(\begin{array}{c}
\rho v \\ \rho v^2 + p(\rho)
\end{array}\right),
\quad
R(U)=\left(\begin{array}{c}
0 \\ \rho v
\end{array}\right)
$$
and
$
Q=(1~0). 
$
In this case, the local equilibria $u=\rho$ are scalar-valued, with 
$$
\Ecal(u) = \left(\begin{array}{c}
\rho \\ 0
\end{array}\right), 
$$
and we do have $QF(\Ecal(u))=0$.

The diffusive regime for the Euler equations with friction is analyzed as follows.
First, according to the general theory, equilibrium solutions satisfy 
$$
\del_t\rho = -\del_x\Big( Q A(\Ecal(u)) \, U_1 \Big),
$$
where
$$
D_U F(\Ecal(u)) = \left(\begin{array}{cc}
0 & 1 \\ p'(\rho) & 0
\end{array}\right).
$$
Here, $U_1$ is the unique solution to
$$
\begin{aligned}
B (\Ecal(u)) U_1 & = -\del_x \big( F(\Ecal(u)) \big),
\\
QU_1 & = 0
\end{aligned}
$$
with 
$$
B (\Ecal(u)) = \left(\begin{array}{cc}
0 & 0 \\ 0 & 1
\end{array}\right), 
\quad\quad
\del_x \Big( F(\Ecal(u)) \Big) = \left(\begin{array}{c}
0 \\ \del_x \big( p(\rho)\big) 
\end{array}\right). 
$$
This leads us to the effective diffusion equation for the Euler equations with friction:  
\be
\label{Euler2}
\del_t\rho = \del^2_x \big( p(\rho) \big), 
\ee
which is a nonlinear parabolic equation, at least away from vacuum, since $p'(\rho) >0$ by assumption. 
Interestingly, at vacuum, this equation may be degenerate since  $p'(\rho)$ typically vanishes at $\rho=0$. 
For instance, in the case of polytropic gases $p(\rho) = \kappa \rho^\gamma$ with 
$\kappa>0$ and $\gamma \in (1,\gamma)$ we obtain 
\be
\label{Euler3}
\del_t\rho = \kappa \gamma \, \del_x \big( \rho^{\gamma-1} \del \rho \big). 
\ee

In addition, by defining the nternal energy $e(\rho)>0$ by 
$$
e'(\rho)=\frac{p(\rho)}{\rho^2},
$$
we easily check that all smooth solutions to \eqref{Euler1} satisfy
\be
\label{Euler3} 
\eps \, \del_t \Big(\rho\frac{v^2}{2}+\rho e(\rho)\Big)
+
\del_x\Big(\rho\frac{v^3}{2}+ (\rho e(\rho) + p(\rho)) v \Big)
=-\frac{\rho v^2}{\eps}, 
\ee
so that the function 
$$
\Phi(U)=\rho\frac{v^2}{2}+\rho e(\rho)
$$
 is convex entropy compatible with the relaxation. 
All the conditions of the general framework are therefore satisfied by the Euler equations with friction.


\subsection{$M1$ model of radiative transfer}

Our next model of interest arises in the theory of radiative transfer, i.e. 
\be
\label{M11} 
\begin{aligned}
& \eps \, \del_t e +\del_x f =\frac{\tau^4-e}{\eps},
\\
& \eps \, \del_t f+\del_x \Big(  \chi\left(f/e\right)e \Big) =
-\frac{f}{\eps},
\\
& \eps \, \del_t\tau=\frac{e-\tau^4}{\eps},
\end{aligned}
\ee
where the radiative energy $e>0$ and the radiative flux $f$ are the main unknowns, 
restricted by  the condition $| f/e | \le 1$, 
while $\tau>0$ denotes the temperature. 
The function $\chi:[-1,1]\to\RR^+$ is called the Eddington factor and, typically,  
$$
\chi(\xi) =\frac{ 3+4\xi^2}{5+2\sqrt{4-3\xi^2}}.
$$
 
This system fits into our late-time/stiff relaxation framework if we introduce 
$$
U=\left(\begin{array}{c}
e \\ f \\ \tau
\end{array}\right),
\qquad
F(U)=\left(\begin{array}{c}
f \\ \chi(\frac{f}{e})e \\ 0
\end{array}\right), 
\quad\quad
R(U)=\left(\begin{array}{c}
e-\tau^4 \\ f \\ \tau^4-e
\end{array}\right).
$$
Now, the equilibria are given by $u=\tau+\tau^4$, with 
$$
\Ecal(u)=\left(\begin{array}{c}
\tau^4 \\ 0 \\ \tau
\end{array}\right),
\qquad \quad
Q := (1~0~1). 
$$ 
We do have $QF(\Ecal(u))=0$ and the assumptions in Sections~2 are satisfied.  

We determine the diffusive regime for the $M1$ model from the expression  
$$
\big( D_U F \big) (\Ecal(u)) = 
\left(\begin{array}{ccc}
0 & 1 & 0 \\ 
\chi(0) & \chi'(0) & 0 \\
0 & 0 & 0
\end{array}\right)
=
\left(\begin{array}{ccc}
0 & 1 & 0 \\ 
\frac{1}{3} & 0 & 0 \\
0 & 0 & 0
\end{array}\right), 
$$
where $U_1$ is the solution to the linear problem 
$$
\begin{aligned}
\left(\begin{array}{ccc}
1 & 0 & -4\tau^3 \\ 
0 & 1 & 0 \\
-1 & 0 & 4\tau^3
\end{array}\right)U_1 & =
\left(\begin{array}{c}
0 \\ 
\del_x\big( \tau^4/3 \big) \\
0 
\end{array}\right),\\
 (1~0~1) U_1 & =0. 
\end{aligned}
$$
Therefore, we have  
$$
U_1 = \left(\begin{array}{c}
0 \\ 
\frac{4}{3} \tau^3\del_x\tau \\
0
\end{array}\right), 
$$
and the effective diffusion equation for the $M1$ system reads  
\be
\label{M12} 
\del_t(\tau+\tau^4) = \del_x\left(
\frac{4}{3}\tau^3\del_x \tau \right).
\ee
Again, an entropy can be associated to this model. 


\subsection{Coupled Euler/$M1$ model}

By combining the previous two examples we arrive at a more involved model: 
\be
\label{EM1} 
\begin{aligned}
& \eps\del_t\rho +\del_x \big( \rho v \big) = 0,
\\
& \eps\del_t\rho v+\del_x \big( \rho v^2+p(\rho) \big) 
= -\frac{\kappa}{\eps} \rho v+\frac{\sigma}{\eps} f,
\\
& \eps\del_t e +\del_x f =0,
\\
& \eps\del_t f+\del_x \Big(  \chi\left(\frac{f}{e}\right)e\Big)
 = -\frac{\sigma}{\eps}f,
\end{aligned}
\ee
in which the same notation as before is used and $\kappa$ and $\sigma$ are positive constants. In the applications, a typical choice for the pressure
is 
$$
p(\rho)=C_p \rho^\eta, \qquad 
C_p\ll 1, 
\qquad 
\eta>1.
$$
To fit this model within the late-time/stiff relaxation framework, we need to set 
$$
U=\begin{pmatrix}\rho\\ \rho v\\ e\\f \end{pmatrix},
\qquad
F(U)=\begin{pmatrix}\rho v\\ \rho v^2+p(\rho)\\ f\\ \chi(\frac{f}{e})e \end{pmatrix},
\qquad
R(U)=\begin{pmatrix}0\\ \kappa \rho v-\sigma f\\ 0\\ \sigma f \end{pmatrix}.
$$
The local equilibria are given by 
$$
 \Ecal(u)=\begin{pmatrix}\rho\\ 0\\ e\\0 \end{pmatrix},
 \qquad
u=QU=\begin{pmatrix}\rho\\ e\end{pmatrix},
\qquad
Q=\begin{pmatrix}1&0&0&0\\0&0&1&0\end{pmatrix}, 
$$
and once again, one has$QF(\Ecal(u))=0$. 

We can then compute  
$$
D_U F(\Ecal(u))=\begin{pmatrix}0&1&0&0\\p'(\rho)&0&0&0\\0&0&0&1\\0&0&\frac{1}{3}&0\end{pmatrix},
\qquad
U_1=\begin{pmatrix}0\\ \frac{1}{\kappa}\Bigl(-\del_x p(\rho)-\frac{1}{3}\del_x e\Bigr)\\0
\\-\frac{1}{3\sigma}\del_x e\end{pmatrix}, 
$$ 
and we arrive at the effective diffusion system for the coupled Euler/$M1$ model: 
\be
\label{EM2} 
\begin{aligned}
&\del_t \rho-\frac{1}{\kappa}\del_x^2 p(\rho)-\frac{1}{3\kappa}\del_x^2 e=0,\\
&\del_t e-\frac{1}{3\sigma}\del_x^2 e=0.
\end{aligned}
\ee
The second equation is a standard heat equation, and its solution serves as a source-term in the first equation. 

 
\subsection{Shallow water with nonlinear friction}

Our final examples requires the more general theory of nonlinear relaxation mentioned in Remark~\ref{rmk26}, and reads 
\be
\label{EN1} 
\begin{aligned}
& \eps \del_t h +\del_x \big( hv \big) 
=0,\\
& \eps \del_t \big( hv \big) 
+\del_x \Big( h \, v^2 + p(h) \Big) = -\frac{\kappa^2(h)}{\eps^2} \, g \, hv|hv|,
\end{aligned}
\ee
where $h$ denotes the fluid height and $v$ the fluid velocity $v$. The pressure is taken to be $p(h)=g\, h^2/2$. 
and $g>0$ is called the gravity constant. The friction coefficient $\kappa:\RR^+ \to\RR^+$ is a positive function, and 
a standard choice is $\kappa(h)=\frac{\kappa_0}{h}$ with $\kappa_0>0$.

The nonlinear version of the late-time/stiff relaxation framework applies if we set 
$$
U=\left(\begin{array}{c}
h \\ hv
\end{array}\right),
\quad
F(U)=\left(\begin{array}{c}
hv \\ hv^2+p(h)
\end{array}\right), 
\quad
R(U)=\left(\begin{array}{c}
0 \\ \kappa^2(h)ghv|hv|
\end{array}\right).
$$
The scalar equilibria $u=h$ are associated with 
$$
\Ecal(u)=\left(\begin{array}{c}
h \\ 0
\end{array}\right), 
\qquad \qquad 
Q=(1~0). 
$$ 
Here, the relaxation is nonlinear and satisfies 
$$
R(\Ecal(u)+\eps U) = \eps^2 R\big( \Ecal(U) + M(\eps)U \big),
$$
with 
$$
M(\eps) := \left(\begin{array}{cc}
\eps & 0 \\ 0 & 1
\end{array}\right).
$$ 
in turn , we may derive a {\sl nonlinear} effective  equation for the Euler equations with nonlinear friction, that is, 
\be
\label{EN2} 
\del_t h = \del_x\left( {\sqrt{h} \over \kappa(h)} \, 
\frac{\del_x h}{\sqrt{|\del_x h|}}\right),
\ee
which is a nonlinear parabolic equation. 

In addition, by introducing the internal energy $e(h) :=gh/2$, we observe that all smooth solutions to \eqref{EN1} 
satisfy the entropy inequality 
\be
\label{EN2} 
\eps\del_t\left( h\frac{v^2}{2}+g\frac{h^2}{2} \right)
+
\del_x\left(
h\frac{v^2}{2}+gh^2
\right)v=-\frac{\kappa^2(h)}{\eps^2}ghv^2|hv|.
\ee
The entropy 
$$
\Phi(U) := h\frac{v^2}{2}+g\frac{h^2}{2}
$$ 
satisfies the compatibility properties relevant to the nonlinear late-time/stiff relaxation theory, with in particular 
$$
R(\Ecal(u)+M(0)\bar{U}_1) = \left(\begin{array}{c}
0 \\ \del_x p(h)
\end{array}\right),
$$
where $\bar{U}_1=(0~\beta)\, $.   We obtain here 
$R(\Ecal(u)+M(0)\bar{U}_1)=c(u)\bar{U}_1$ with
$$
c(u) = g\kappa(h)\sqrt{h|\del_x h|}\ge 0.
$$


\section{Asymptotic-preserving finite volume schemes} 

\subsection{General strategy}

We are going now to design finite volume schemes that are consistent with the asymptotic regime $\eps \to 0$
determined in the previous section and, indeed, to recover an effective diffusion equation
that is independent of the  mesh-size. The discretization of hyperbolic-to-hyperbolic regimes were 
investigated first by Ji and Xin \cite{JX}. Here we propose a framework to cover hyperbolic-to-parabolic regimes. 
For earlier work on this latter regime see \cite{BCD,BT,BOP,BC,BD}.

\begin{itemize}

\item[Step 1.]  Our construction is based on a standard finite volume scheme for  the homogeneous system
$$ 
\del_t U+\del_x F(U) =0, 
$$
and we will begin by describing such a scheme. 

\item[Step 2.] We then will modify the above scheme and include a matrix-valued free parameter, 	allowing us to 
approximate the non-homogeneous system 
$$
\del_t U+\del_x F(U) = -\gamma \, R(U),
$$
for a fix coefficient $\gamma>0$. 

\item[Step 3.] Finally, we will perform an asymptotic analysis after replacing the discretization parameter
 $\Delta t$ by $\eps \Delta t$, and 
$\gamma$ by $1/\eps$. Our goal then will be to determine the free parameters so that to ensure the asymptotic-preserving property.
\end{itemize}

Let us briefly define the so-called HLL discretization of the homogeneous system, as proposed by 
Harten, Lax, and van Leer \cite{HLL}. For simplicity we present here the solver based on 
a single constant intermediate state. The mesh is assumed to be uniform mesh made of cells of length $\Delta x$: 
$$
[x_{i-\demi},x_{i+\demi}], 
\qquad 
x_{i+\demi}=x_i+\frac{\Delta x}{2} 
$$
for all $i=\ldots, -1,0,1, \ldots$.  
The time discretization is based on a parameter $\Delta t$, restricted by the famous CFL condition (Courant, Friedrichs, Lewy, cf.~\cite{CFL}) with  
$$
t^{m+1}=t^m+\Delta t. 
$$

Starting from som initial data (lying in the convex set $\Omega$): 
$$
U^0(x) = \frac{1}{\Delta x}
\int_{x_{i-\demi}}^{x_{i+\demi}} U(x,0)dx,
\qquad x\in[x_{i-\demi},x_{i+\demi}).
$$
we construct piecewise constant approximations at each time $t^m$: 
$$
U^m(x)=U_i^m, 
\qquad x\in[x_{i-\demi},x_{i+\demi}), \quad i\in\ZZ. 
$$

Following Harten, Lax and van Leer \cite{HLL},  at each cell interface we use the
{\sl approximate Riemann solver:} 
$$
\Ubar_\Rcal(\frac{x}{t};U_L,U_R) =\left\{
\begin{aligned}
& U_L, \quad \frac{x}{t}<-b,\\
& \Ubar^\star, \quad -b<\frac{x}{t}<b,\\
& U_R, \quad \frac{x}{t}>b,
\end{aligned}\right.
$$
where $b>0$ is sufficiently large.
The ``numerical cone'' (and numerical diffusion) is thus determined by the 
parameter $b>0$ and here, for simplicity in the presentation, we have assumed
 a single constant $b$ but, more generally,  one can introduce two speeds $b^-_{i+\demi}<b^+_{i+\demi}$ at each interface. 

We introduce the intermediate state 
$$
\Ubar^\star = \frac{1}{2}(U_L+U_R)-\frac{1}{2b}
\big( F(U_R) - F(U_L) \big),
$$
and assume the CFL condition 
$$
b\frac{\Delta t}{\Delta x}\le\frac{1}{2},
$$
so that the underlying approximate Riemann solutions are non-interacting. Our global approximate solutions 
$$
\Ubar^m_{\Delta x}(x,t^m+t), \qquad t\in[0,\Delta t), \quad x \in \RR. 
$$
are then obtained as follows. 

The approximations at the next time level $t^{m+1}$ are determined from  
$$
\Ubar_i^{m+1} = \frac{1}{\Delta x}
\int_{x_{i-\demi}}^{x_{i+\demi}} \Ubar^m_{\Delta x}(x,t^m+\Delta
t) dx.
$$
Then, recalling 
$
\Ubar^\star_{i+\demi} = \frac{1}{2}(U_i^m+U_{i+1}^m)-\frac{1}{2b}
(F(U_{i+1}^m)-F(U_i^m)),
$
and integrating out the expression given by the Riemann solutions, 
we arrive at the {\sl scheme for the homogeneous system}
$$
\Ubar_i^{m+1} = U_i^m - \frac{\Delta t}{\Delta x}
\Big( F^{HLL}_{i+\demi} - F^{HLL}_{i-\demi} \Big),
$$
where
$$
F^{HLL}_{i+\demi} = \frac{1}{2} \Big( F(U_i^m)+F(U_{i+1}^m) \Big)
-\frac{b}{2}(U_{i+1}^m-U_i^m). 
$$
(More generally one could take into account two speeds $b^-_{i+\demi}<b^+_{i+\demi}$.)

We observe that the above scheme enjoys invariant domains.
The intermediate states $\Ubar^\star_{i+\demi}$ can be written in the form of a convex combination 
$$
\Ubar^\star_{i+\demi}= 
\frac{1}{2}\left(U_i^m+\frac{1}{b}F(U_i^m)\right)
+ 
\frac{1}{2}\left(U_{i+1}^m-\frac{1}{b}F(U_{i+1}^m)\right) \in \Omega,
$$
provided $b$ is large enough. An alternative decomposition is given by  
$$
\Ubar^\star_{i+\demi}
= 
\frac{1}{2}\Big( I + \frac{1}{b} \Abar(U_i^m, U^m_{i+1}) \Big) \, U_i^m
+ 
\frac{1}{2}\Big( I - \frac{1}{b} \Abar(U_i^m, U^m_{i+1}) \Big) \, U_{i+1}^m, 
$$
where $\Abar$ is an  ``average'' of $D_U F$.   
By induction, it follows that  $\Ubar_i^m$ in $\Omega$ for all $m,i$.

\subsection{Discretization of the relaxation term} 

We start from the following {\sl modified Riemann solver:} 
$$
U_\Rcal(\frac{x}{t};U_L,U_R) =\left\{
\begin{aligned}
& U_L, \quad \frac{x}{t}<-b,\\
& U^{\star L}, \quad -b<\frac{x}{t}<0,\\
& U^{\star R}, \quad 0<\frac{x}{t}<b,\\
& U_R, \quad \frac{x}{t}>b,
\end{aligned}\right.
$$
with the following states at the interface: 
$$
\begin{aligned}
& U^{\star L} = \alphab\Ubar^\star +  (I-\alphab)\big( U_L - \bar{R}(U_L) \big),
\\
& U^{\star R} = \alphab\Ubar^\star + 
    (I-\alphab) \big( U_R - \bar{R}(U_R) \big).
\end{aligned}
$$
We have here introduced some $N\times N$--matrix  and $N$--vector defined by 
$$
\alphab = \left( I +\frac{\gamma\Delta
  x}{2b}(I+\underline{\sigma})
\right)^{-1},
\qquad  
\bar{R}(U) = (I+\underline{\sigma})^{-1} R(U). 
$$
The term $\underline{\sigma}$ is a parameter matrix to be chosen 
so that (all inverse matrices are well-defined and) the correct asymptotic regime is recovered
at the discrte level.

At each interface $x_{i+\demi}$, we use the modified Riemann
solver 
$$
U_{\Rcal}(\frac{x-x_{i+\demi}}{t-t^m}; U_i^m,U_{i+1}^m)
$$ 
and we superimpose non-interacting Riemann solutions 
$$
U^m_{\Delta x}(x,t^m+t), \qquad t\in[0,\Delta t), \quad x \in \RR.
$$ 
The new approximate solution at the next time $t^{m+1}$ is 
$$
U_i^{m+1} = 
\int_{x_{i-\demi}}^{x_{i+\demi}} U^m_{\Delta x}(x,t^m+\Delta t)dx.
$$
By integrating out the Riemann solutions, we arrive at the 
discretized balance law
\be
\label{904} 
\begin{aligned}
& \frac{1}{\Delta t}(U_i^{m+1}-U_i^m) +\frac{1}{\Delta x}
\Big(\alphab_{i+\demi} F^{HLL}_{i+\demi}
-\alphab_{i-\demi} F^{HLL}_{i-\demi} \Big)
\\
& =
\frac{1}{\Delta x}
  (\alphab_{i+\demi}-\alphab_{i-\demi}) F(U_i^m)
-\frac{b}{\Delta x} (I-\alphab_{i-\demi})
  \bar{R}_{i-\demi}(U_i^m)
\\
& \quad - 
\frac{b}{\Delta x} (I-\alphab_{i+\demi})
  \bar{R}_{i+\demi}(U_i^m). 
\end{aligned}
\ee

Observing that the discretized source can rewritten as 
$$
\aligned
\frac{b}{\Delta x}(I-\alphab_{i+\demi})
  \bar{R}_{i+\demi}(U_i^m) 
& = 
\frac{b}{\Delta x}\alphab_{i+\demi}
 (\alphab_{i+\demi}^{-1}-I)
  \bar{R}_{i+\demi}(U_i^m)
\\
& = \frac{\gamma}{2}\alphab_{i+\demi}R(U_i^m) 
\endaligned
$$
and, similarly, 
$$
\frac{b}{\Delta x}(I-\alphab_{i-\demi})
  \bar{R}_{i-\demi}(U_i^m) = 
\frac{\gamma}{2}\alphab_{i-\demi}R(U_i^m),
$$
we conclude that the proposed
 {\sl finite volume schemes for late-time/stiff-relaxation problems} takes the form 
\be
\label{914} 
\begin{aligned}
& \frac{1}{\Delta t}(U_i^{m+1}-U_i^m) + \frac{1}{\Delta x}
(\alphab_{i+\demi} F^{HLL}_{i+\demi}
-\alphab_{i-\demi} F^{HLL}_{i-\demi})
\\
& =
\frac{1}{\Delta x}
  (\alphab_{i+\demi}-\alphab_{i-\demi}) F(U_i^m)
-\frac{\gamma}{2}
(\alphab_{i+\demi}+\alphab_{i-\demi})R(U_i^m).
\end{aligned}
\ee

\begin{theorem}[Properties of the finite volume scheme]
Provided 
$$
\underline{\sigma}_{i+\demi}-\underline{\sigma}_{i-\demi}=\Ocal(\Delta
x)
$$ 
and the matrix-valued map $\underline{\sigma}$ is sufficiently smooth, 
then the modified finite volume scheme is {\rm consistent} with the hyperbolic system with relaxation.
The following {\rm invariant domain} property holds:  
provided all the states at the interfaces 
$$
\begin{aligned}
&  U^{\star L}_{i+\demi} = \alphab_{i+\demi}\Ubar^\star_{i+\demi} + 
    (I-\alphab_{i+\demi})(U_i^m-\bar{R}(U_i^m)),\\
& U^{\star R}_{i+\demi} = \alphab_{i+\demi}\Ubar^\star_{i+\demi} + 
    (I-\alphab_{i+\demi})(U_{i+1}^m-\bar{R}(U_{i+1}^m))
\end{aligned}
$$
belong to $\Omega$, then all the states $U_i^m$ belong to $\Omega$.
\end{theorem}

\subsection{Discrete late-time asymptotic regime}

As explained earlier, we now replace  $\Delta t$ by $\Delta t/\eps$ and $\gamma$ by $1 / \eps$, and consider 
$$
\begin{aligned}
& \frac{\eps}{\Delta t}(U_i^{m+1}-U_i^m) +\frac{1}{\Delta x}
(\alphab_{i+\demi} F^{HLL}_{i+\demi}
-\alphab_{i-\demi} F^{HLL}_{i-\demi})
\\
&
= \frac{1}{\Delta x}
  (\alphab_{i+\demi}-\alphab_{i-\demi}) F(U_i^m)
-\frac{1}{2\eps}
(\alphab_{i+\demi}+\alphab_{i-\demi})R(U_i^m),
\end{aligned}
$$
where
$$
\alphab_{i+\demi}=\left(
I+\frac{\Delta x}{2\eps b}(I+\underline{\sigma}_{i+\demi})
\right)^{-1}.
$$
Plugging in an expansion near an equilibrium state 
$$
U_i^m = \Ecal(u_i^m) + \eps(U_1)_i^m + \Ocal(\eps^2),
$$
we find 
$$
\aligned 
F^{HLL}_{i+\demi}  
& =
{1 \over 2} F\big( \Ecal(u_i^m) \big) + {1 \over 2} F\big( \Ecal(u_{i+1}^m) \big) 
-\frac{b}{2}\left( \Ecal(u_{i+1}^m) - \Ecal(u_i^m) \right)
+ \Ocal(\eps),
\\
\frac{1}{\eps}R(U_i^m) & = B(\Ecal(u_i^m)) (U_1)_i^m +
  \Ocal(\eps), 
\\
\alphab_{i+\demi}
& =
\frac{2b \eps}{\Delta x} \big(I+\underline{\sigma}_{i+\demi} \big)^{-1} + \Ocal(1).
\endaligned
$$

The {\sl first-order} terms for the discrete scheme leads us to  
$$
\begin{aligned}
& \frac{1}{\Delta t}(\Ecal(u_i^{m+1}) - \Ecal(u_i^m))
\\
& = - \frac{2b}{\Delta x^2}\left(
(I+\underline{\sigma}_{i+\demi})^{-1} F^{HLL}_{i+\demi}|_{\Ecal(u)}
-
(I+\underline{\sigma}_{i-\demi})^{-1} F^{HLL}_{i-\demi}|_{\Ecal(u)}
\right)
\\
& \quad + 
\frac{2b}{\Delta x^2} \left(
(I+\underline{\sigma}_{i+\demi})^{-1}
-
(I+\underline{\sigma}_{i-\demi})^{-1}
\right) F(\Ecal(u_i^m))
\\
&
\quad - \frac{b}{\Delta x} \left(
(I+\underline{\sigma}_{i+\demi})^{-1}
+
(I+\underline{\sigma}_{i-\demi})^{-1}
\right)B(\Ecal(u_i^m))(U_1)_i^m.
\end{aligned}
$$
At this juncture, we assume the existence of an $n\times n$ matrix
$\Mcal_{i+\demi}$ such that
$$
Q \big( I+\underline{\sigma}_{i+\demi} \big)^{-1}=\frac{1}{b^2} \, \Mcal_{i+\demi} Q.
$$
We then multiply the equation above by the  $n \times N$ matrix $Q$ and obtain 
$$
\frac{1}{\Delta t}(u_i^{m+1}-u_i^m) =
-\frac{2}{b\Delta x^2}\left(
\Mcal_{i+\demi} Q F^{HLL}_{i+\demi}|_{\Ecal(u)} - 
\Mcal_{i-\demi} Q F^{HLL}_{i-\demi}|_{\Ecal(u)} \right),
$$
with 
$$
\begin{aligned}
Q F^{HLL}_{i+\demi}|_{\Ecal(u)}&=  {Q \over 2}  F(\Ecal(u_i^m)) 
+ {Q \over 2}  F(\Ecal(u_{i+1}^m))
-\frac{b}{2} Q \left(\Ecal(u_{i+1}^m)-\Ecal(u_i^m) \right)
\\
&= -\frac{b}{2} (u_{i+1}^m -u_i^m).
\end{aligned}
$$
The discrete asymptotic system is thus
\be
\label{944} 
\frac{1}{\Delta t}(u_i^{m+1}-u_i^m) =
\frac{1}{\Delta x^2}
\Big(
\Mcal_{i+\demi}(u_{i+1}^m -u_i^m) + 
\Mcal_{i-\demi}(u_{i-1}^m -u_i^m)
\Big).
\ee

Recall that for some $n\times n$ matrix $\Mcal(u)$ the effective equation reads
$$
\del_t u = \del_x \left(\Mcal(u) \del_x u \right).
$$

\begin{theorem} [Discrete  late-time asymptotic-preserving property] 
Assume the following conditions on the matrix-valued coefficients: 

\begin{itemize} 

\item[$\bullet$] The matrices $$
I+\underline{\sigma}_{i+\demi}, \qquad \Big(1+\frac{\Delta x}{2\eps b} \Big)I+\underline{\sigma}_{i+\demi}
$$
are 
invertible for $\eps \in [0,1]$.

\item[$\bullet$] There exists an $n\times n$ matrix $\Mcal_{i+\demi}$ ensuring the commutation condition 
$$
Q(I+\underline{\sigma}_{i+\demi})^{-1}=\frac{1}{b^2}\Mcal_{i+\demi} Q. 
$$
\item[$\bullet$] The discrete form of $\Mcal(u)$ at each cell interface $x_{i+\demi}$ satisfies 
$$
\Mcal_{i+\demi} = \Mcal(u) + \Ocal(\Delta x). 
$$
\end{itemize}
Then, the effective system associated with the discrete scheme coincides with the one
of the late-time/stiff relaxation framework. 
\end{theorem}

 We refer the reader to \cite{BLT} for numerical experiments with this scheme, which turns out to 
efficiently compute the late-time behavior of solutions. It is observed therein that 
asymptotic solutions may have large gradients but are actually regular. 
We also emphasize that we rely here on the CFL stability condition based on the homogeneous hyperbolic system, i.e., a restriction 
on $\Delta t/\Delta x$ only is imposed.
In typical tests, about $10 000$ time-steps were used to reach the late-time behavior
and,  for simplicity, the initial data were taken in the image of $Q$. A {\sl reference} solution, needed for a comparison, 
was obtained by solving the parabolic equation, under a (stronger) restriction on $\Delta t/(\Delta x)^2$. 

The proposed theoretical framework for late-time/stiff relaxation problems thus led us 
to the development of a good strategy
to design asymptotic-preserving schemes involving matrix-valued  parameter. 
The convergence analysis ($\eps \to 0$) and the numerical analysis ($\Delta x \to 0$) for the problems under consideration 
are important and challenging open problems.  
It will also very interesting to apply our technique to, for instance, 
plasma mixtures in a  multi-dimensional setting.  

\section*{Acknowledgments}
 
The author was partially supported by  
the Agence Nationale de la Recherche (ANR) through the grant 06-2-134423, and 
by the Centre National de la Recherche Scientifique (CNRS).


\end{document}